\documentclass[twoside,draft,12pt]{article}
\usepackage{amssymb}

\usepackage{makeidx}
\usepackage{a4wide}
\usepackage{amsthm}
\usepackage{amsmath}
\usepackage[all]{xy}
\usepackage{enumerate}
\usepackage{fancyhdr}

\newtheorem{theorem}{Theorem}[section]
\newtheorem{proposition}[theorem]{Proposition}
\newtheorem{corollary}[theorem]{Corollary}
\newtheorem{lemma}[theorem]{Lemma}
\theoremstyle{definition}

\newtheorem{c-example}[theorem]{Counter Example}

\newtheorem*{Beweis}{Proof}
\newtheorem{definition}[theorem]{Definition}
\newtheorem{punto}[theorem]{}
\theoremstyle{remark}
\newtheorem{remark}[theorem]{Remark}

\CompileMatrices
\input{tcilatex}

\begin{document}

\title{On the Linear Weak Topology and Dual Pairings over Rings\thanks{%
MSC (2000): 13Jxx, 16E60, 16E99, 16W80 \newline
Keywords: Finite Topology, Linear Weak Topology, Dual Pairings, Dense
Pairings, Locally Projective Modules}}
\author{\textbf{Jawad Y. Abuhlail} \\
Department of Mathematical Sciences\\
King Fahd University of Petroleum $\&$ Minerals\\
31261 Dhahran - Saudi Arabia\\
abuhlail@kfupm.edu.sa}
\date{}
\maketitle

\begin{abstract}
In this note we study the weak topology on paired modules over a (not
necessarily commutative) ground ring. Over QF rings we are able to recover
most of the well known properties of this topology in the case of
commutative base fields. The properties of the linear weak topology and the
dense pairings are then used to characterize pairings satisfying the so
called $\alpha $-condition.
\end{abstract}

\section*{Introduction}

Let $R$ be a commutative field, $V,W$ be vector spaces over $R$ with a
non-degenerating $R$-bilinear form $\beta :V\times W\rightarrow R,$ $P:=(V,W)
$ be the induced $R$-pairing and consider $V\overset{\kappa _{P}}{%
\hookrightarrow }W^{\ast }$ and $W\overset{\chi _{P}}{\hookrightarrow }%
V^{\ast }$ as vector subspaces. For every subset $K\subseteq W$
(respectively $X\subseteq W^{\ast }$) set 
\begin{equation*}
\mathrm{An}(K)=\{f\in W^{\ast }\mid f(K)=0\}\text{ (respectively }\mathrm{Ke}%
(X)=\bigcap \{\mathrm{Ker}(f)\mid f\in X\}\text{).}
\end{equation*}
Considering $R$ with the discrete topology, $R^{W}$ with the product
topology, the induced relative topology on $V\subseteq W^{\ast }\subset R^{W}
$ is called the \emph{linear weak topology} $V[\frak{T}_{ls}(W)]$ and has
basis of neighbourhoods of $0_{V}:$%
\begin{equation*}
\{F^{\bot }:V\cap \mathrm{An}(F)\mid F=\{w_{1},...,w_{k}\}\subset W\text{ a
finite subset}\}.
\end{equation*}
The linear weak topology on $W\subseteq V^{\ast }\subset R^{V}$ is defined
analogously. The closure of any vector subspace $X\subset V$ is given by $%
\overline{X}:=X^{\bot \bot }.$ A closed (open) vector subspace $X\subset V$
has the form $X=K^{\bot },$ where $K\subset W$ is any (finite dimensional)
vector subspace. The embeddings $W\hookrightarrow V^{\ast }$ respectively $%
V\hookrightarrow W^{\ast }$ imply that $V\subseteq W^{\ast }$ respectively $%
W\subseteq V^{\ast }$ are dense. The properties of this topology are well
known and were studied by several authors (e.g. \cite{Kot66}, \cite{KN76}, 
\cite{Rad73}, \cite{LR97}).

For the case of arbitrary base rings most of the properties of this topology
(including the characterizations of closure, closed, open and dense
submodules) are not valid anymore. The aim of this note is to study the
properties of this topology induced on paired modules over arbitrary ground
rings. In particular we extend results obtained by the author \cite[Anhang]
{Abu2001} on this topology from the case of commutative base rings to the
arbitrary case. In contrast to the proofs in the case of base fields, which
depend heavily on the existence of bases, our proofs are in module theoretic
terms.

Throughout this note $R$ denotes a (\emph{not necessarily} \emph{commutative}%
)\emph{\ }associative ring with $1_{R}\neq 0_{R}.$ We consider $R$ as a
right (and a left) linear topological ring with the discrete topology. The
category of unitary left (right) $R$-modules will be denoted with $_{R}%
\mathcal{M}$ ($\mathcal{M}_{R}$). The category of unitary $R$-bimodules is
denoted with $_{R}\mathcal{M}_{R}.$ For a right (a left) $R$-module $L$ we
denote with $L^{\ast }$ ($^{\ast }L$) the set of all $R$-linear mappings
from $L$ to $R.$ If $V$ is an $R$-module, then an $R$-submodule $X\subset V$
is called $R$-cofinite, if $V/X$ is finitely generated as an $R$-module.

Let $L$ be a right\ (a left) $R$-module and $K\subset L$ be an $R$%
-submodule. We call $K\subset L$ $N$-pure for some left (right) $R$-module $%
N,$ if the canonical mapping $\iota _{K}\otimes id_{N}:K\otimes
_{R}N\rightarrow L\otimes _{R}N$ ($id_{N}\otimes \iota _{L}:N\otimes
_{R}K\rightarrow N\otimes _{R}L$) is an embedding. We call $K\subset L$ pure
(in the sense of Cohn), if $K\subset L$ is $N$-pure for every left (right) $%
R $-module $N.$

\section{The linear weak topology}

\begin{punto}
$R$\textbf{-pairings. }A \emph{left }$R$\emph{-pairing} $P=(V,W)$ consists
of a left $R$-module $W$ and a right $R$-module $V$ with an $R$-linear
mapping $\kappa _{P}:V\rightarrow $ $^{\ast }W$ (equivalently $\chi
_{P}:W\rightarrow V^{\ast }$). For left $R$-pairings $(V,W),$ $(V^{\prime
},W^{\prime })$ a morphism $(\xi ,\theta ):(V^{\prime },W^{\prime
})\rightarrow (V,W)$ consists of $R$-linear mappings $\xi :V\rightarrow
V^{\prime }$ and $\theta :W^{\prime }\rightarrow W,$ such that 
\begin{equation}
<\xi (v),w^{\prime }>=<v,\theta (w^{\prime })>\text{ for all }v\in V\text{
and }w^{\prime }\in W^{\prime }.  \label{comp}
\end{equation}
Let $P=(V,W)$ be a left $R$-pairing, $V^{\prime }\subset V$ be a right $R$%
-submodule, $W^{\prime }\subset W$ be a (pure) left $R$-submodule with $%
<V^{\prime },W^{\prime }>=0.$ Then $Q:=(V/V^{\prime },W^{\prime })$ is a
left $R$-pairing, $(\pi ,\iota ):(V/V^{\prime },W^{\prime })\rightarrow (V,W)
$ is a morphism of left $R$-pairings and we call $Q\subset P$ a (\emph{pure}%
) \emph{left }$R$\emph{-subpairing.} The left $R$-pairings with the
morphisms defined above build a category, which we denote by $\mathcal{P}%
_{l}.$ The category of right $R$-pairings $\mathcal{P}_{r}$ is defined
analogously.
\end{punto}

\begin{punto}
\textbf{The finite topology.} Consider $R$ with the discrete topology. For
every set $\Lambda $ we consider $R^{\Lambda }$ with the \emph{product
topology} and identify it with the set of all mappings from $\Lambda $ to $R.
$ If $W$ is a left $R$-module, then the induced relative topology on the
right $R$-submodule $^{\ast }W\subset R^{W}$ is called the \emph{finite
topology }and makes $^{\ast }W$ a linear topological right $R$-module with
basis of neighbourhoods of $0_{^{\ast }W}:$%
\begin{equation*}
\mathcal{B}_{f}{\normalsize (0_{^{\ast }W})}:=\{\mathrm{An}(F)|\text{ }%
F=\{w_{1},...,w_{k}\}\subset W\text{ is a finite subset}\}.
\end{equation*}
\end{punto}

\begin{punto}
Let $P=(V,W)$ be a left $R$-pairing and consider the right $R$-submodule $%
^{\ast }W\subset R^{W}$ with the finite topology. Then there is a unique
topology on $V,$ the \emph{linear weak\ topology} $V[\frak{T}_{ls}^{r}(W)],$
such that $\kappa _{P}:V\rightarrow $ $^{\ast }W$ is continuous. A basis of
neighbourhoods of $0_{V}$ is given by the\emph{\ neighbourhoods} 
\begin{equation*}
\mathcal{B}{\normalsize _{f}(0_{V})}:=\{F^{\bot }:=\kappa _{P}^{-1}(\mathrm{%
An}(F))|\text{ }F=\{w_{1},...,w_{k}\}\subset W\text{ is a finite subset}\}.
\end{equation*}
The closure $\overline{X}$ of any subset $X\subseteq V$ is then given by 
\begin{equation*}
\overline{X}=\bigcap \{X+F^{\bot }\mid F\subset W\text{ is a finite subset}%
\}.
\end{equation*}
Analogously one can consider $W$ as a left linear topological $R$-module
with the linear weak topology $W[\frak{T}_{ls}^{l}(V)],$ which is the finest
topology on $W$ that makes $\chi _{W}:W\rightarrow V^{\ast }$ continuous (we
consider $V^{\ast }\subset R^{V}$ with the finite topology).
\end{punto}

\begin{lemma}
Let $P=(V,W)$ be a left $R$-pairing and consider $V$ with the linear weak
topology $V[\frak{T}_{ls}^{r}(W)].$

\begin{enumerate}
\item  $V[\frak{T}_{ls}^{r}(W)]$ is Hausdorff if and only if $V\overset{%
\kappa _{P}}{\hookrightarrow }$ $^{\ast }W.$

\item  If $\kappa _{P}(V)\subset $ $^{\ast }W$ is dense and $_{R}R$ is $W$%
-injective, then $\widehat{V}\simeq $ $^{\ast }W$ \emph{(}where $\widehat{V}$
is the completion of $V$ w.r.t. $V[\frak{T}_{ls}^{r}(W)]$\emph{)}.

\item  The finite topology on $^{\ast }W$ is Hausdorff. If $_{R}R$ is $W$%
-injective, then $^{\ast }W$ is complete.
\end{enumerate}
\end{lemma}

\begin{Beweis}
Denote with $\mathcal{W}^{f}$ the class of all finitely generated $R$%
-submodules of $W.$

\begin{enumerate}
\item  This is evident, while 
\begin{equation*}
\overline{0_{V}}:=\bigcap \{K^{\bot }|\text{ }K\in \mathcal{W}{\normalsize %
^{f}\}=(}\sum \{K\in \mathcal{W}{\normalsize ^{f}\})^{\bot }=W}^{\bot }=%
\mathrm{Ker}{\normalsize (\kappa _{P})}.
\end{equation*}

\item  Consider for every left $R$-submodule $K\overset{\iota _{K}}{%
\hookrightarrow }W$ the $R$-linear mapping 
\begin{equation*}
\varphi _{K}:V\rightarrow \text{ }^{\ast }K,\text{ }v\mapsto \lbrack
k\mapsto <v,k>].
\end{equation*}
Since $_{R}R$ is $W$-injective, $\iota _{K}^{\ast }:$ $^{\ast }W\rightarrow $
$^{\ast }K$ is surjective. By assumption $\kappa _{P}(V)\subset $ $^{\ast }W$
is dense and consequently for every finitely generated left $R$-submodule $%
K\subset W,$ the $R$-linear mapping $\varphi _{K}$ is surjective, hence $%
V/K^{\bot }\simeq $ $^{\ast }K.$ If we write $W=\underrightarrow{lim}%
K_{\lambda }$ as a direct system $\{K_{\lambda }\}_{\Lambda }$ of its
finitely generated $R$-submodules, then 
\begin{equation*}
\widehat{V}:={\underleftarrow{lim}V}/K_{\lambda }^{\bot }\simeq 
\underleftarrow{lim}\text{ }^{\ast }K_{\lambda }\simeq \mathrm{Hom}_{R-}(%
\underrightarrow{lim}K_{\lambda },R)=\text{ }^{\ast }W.
\end{equation*}

\item  The result follows from (1) and (2).$\blacksquare $
\end{enumerate}
\end{Beweis}

\begin{definition}
An $R$-module $U$ is called \emph{FP-injective}, if every diagram of $R$%
-modules 
\begin{equation*}
\xymatrix{ 0 \ar[r]& K \ar[d]_{f} \ar[r] & R^{({\Bbb N})} \ar@{.>}[dl]^{g}
\\ & U & }
\end{equation*}
with exact row and $K$ finitely generated can be completed commutatively
with some $R$-linear mapping $g:R^{(\mathbb{N})}\rightarrow U.$
\end{definition}

An important role in studying the linear weak topology is played by the so
called

\begin{punto}
\textbf{Annihilator conditions} \label{An-Ke}(\cite[28.1]{Wis91}) Let $N$ be
an $R$-module.

\begin{enumerate}
\item  For every $R$-submodule $L\subset N$ we have 
\begin{equation*}
\mathrm{KeAn}(L)=L\Leftrightarrow \text{ }N/L\text{ is }R\text{-cogenerated.}
\end{equation*}

\item  If $R$ is $N$-injective, then 
\begin{equation*}
\mathrm{An}(L_{1}\cap L_{2})=\mathrm{An}(L_{1})+\mathrm{An}(L_{2})\text{ for
all }R\text{-submodules }L_{1},L_{2}\subset N.
\end{equation*}

\item  If $R$ is injective, or if $N$ is finitely generated and $R$ is
FP-injective, then for every finitely generated $R$-submodule $X\subset $ $%
\mathrm{Hom}(N,R)$ we have $\mathrm{AnKe}(X)=X.$
\end{enumerate}
\end{punto}

\qquad We call the ring $R$ a \emph{QF ring}, if $R_{R}$ (equivalently $%
_{R}R $) is Noetherian and a cogenerator (e.g. \cite[48.15]{Wis91}).

\begin{lemma}
\label{orth-clos}Let $P=(V,W)$ be a left $R$-pairing and consider $V$ with
the linear weak topology $V[\frak{T}_{ls}^{r}(W)].$

\begin{enumerate}
\item  $\overline{X}\subseteq X^{\bot \bot }$ for any subset $X\subset V.$
Consequently every orthogonally closed right $R$-submodule of $V$ is closed.

\item  If $R_{R}$ is Noetherian, then all open right $R$-submodules of $V$
are $R$-cofinite.

\item  Let $X\subset V$ be a right $R$-submodule, so that $V/X$ is $R$%
-cogenerated. If $\mathrm{An}(X)=\chi _{P}(X^{\bot }),$ then $X$ is closed.
If moreover $_{R}R$ is Noetherian, $X\subset V$ is $R$-cofinite and $W%
\overset{\chi _{P}}{\hookrightarrow }V^{\ast },$ then $X$ is open.

\item  Let $R_{R}$ be Artinian.

\begin{enumerate}
\item  A right $R$-submodule $X\subset V$ is open if and only if it is
closed and $R$-cofinite.

\item  Let $X\subset Y\subset V$ be right $R$-submodules. If $X\subset V$ is
closed and $R$-cofinite, then $Y\subset V$ is also closed and $R$-cofinite.
\end{enumerate}

\item  Assume $V\subseteq $ $^{\ast }W.$

\begin{enumerate}
\item  If $_{R}R$ is injective, or if $_{R}W$ is finitely generated and $%
_{R}R$ is FP-injective, then every finitely generated right $R$-submodule $%
X\subset V$ is closed.

\item  Let $V_{R}$ be finitely generated. If $_{R}R$ is injective and $R_{R}$
is Noetherian \emph{(}e.g. $R$ is a QF ring\emph{)}, then all right $R$%
-submodules of $V$ are closed.
\end{enumerate}
\end{enumerate}
\end{lemma}

\begin{Beweis}
\begin{enumerate}
\item  Let $\widetilde{x}\in \overline{X}$ be arbitrary. For every $w\in
X^{\bot }$ there exist $x_{w}\in X$ and $v_{w}\in \{w\}^{\bot }$ with $%
\widetilde{x}=x_{w}+v_{w}$ and so $<\widetilde{x},w>=0.$ Consequently $%
\overline{X}\subseteq X^{\bot \bot }.$ If $X$ is orthogonally closed, then $%
\overline{X}\subseteq X^{\bot \bot }=X,$ i.e. $X$ is closed.

\item  Let $X\subset V$ be an open right $R$-submodule. By definition there
exists a finitely generated left $R$-submodule $K\subset W,$ such that $%
K^{\bot }\subset X.$ If $R_{R}$ is Noetherian, then $^{\ast }K_{R}$ is
finitely generated, hence $K^{\bot }\subset V$ is $R$-cofinite. Consequently 
$X\subset V$ is $R$-cofinite.

\item  Let $X\subset V$ be a right $R$-submodule, so that $V/X$ is $R$%
-cogenerated. If $\mathrm{An}(X)=\chi _{P}(X^{\perp }),$ then it follows by 
\ref{An-Ke} (1) that 
\begin{equation*}
X=\mathrm{KeAn}(X)=\mathrm{Ke}(\chi _{P}(X^{\perp }))=X^{\perp \perp }.
\end{equation*}
By (1) $X$\ is closed. Assume now that $_{R}R$ is Noetherian, $X\subset V$
is $R$-cofinite and $W\overset{\chi _{P}}{\hookrightarrow }V^{\ast }.$ Then
by assumption $X^{\bot }=\mathrm{An}(X)\simeq (V/X)^{\ast }$ is finitely
generated in $_{R}\mathcal{M}$ and so $X=(X^{\bot })^{\bot }$ is open.

\item  Assume $R_{R}$ to be Artinian and let $X\subset V$ be a right $R$%
-submodule.

\begin{enumerate}
\item  Every open $R$-submodule $X\subset V$ is closed without any
assumptions on $R$ and is $R$-cofinite by (2). On the other hand, let $%
X\subset V$ be $R$-cofinite and closed. Since $R_{R}$ is Artinian $V/X$ is
finitely cogenerated (e.g. \cite[31.4]{Wis91}), hence open by \cite[1.8]
{Ber94}.

\item  Let $X\subset V$ be $R$-cofinite and closed. Then $X$ is by (a) open
and so $Y\supset X$ is open, hence closed. Obviously $Y\subset V$ is $R$%
-cofinite.
\end{enumerate}

\item  Let $V\overset{\kappa _{P}}{\hookrightarrow }$ $^{\ast }W$ be an
embedding.

\begin{enumerate}
\item  If $X\subset V$ is a finitely generated right $R$-submodule, then we
have under our assumptions and applying \ref{An-Ke} (3): $X^{\bot \bot
}=V\cap \mathrm{AnKe}(X)=X,$ hence $X$ is closed by (1).

\item  Since $V_{R}$ is finitely generated and $R_{R}$ is Noetherian, \emph{%
all} right $R$-submodules of $V$ are finitely generated. Since, by
assumption, $_{R}R$ is injective, the result follows by (a).$\blacksquare $
\end{enumerate}
\end{enumerate}
\end{Beweis}

\section*{Closed and open submodules}

For a left $R$-pairing $(V,W)$ we characterize in what follows the closed
(the open) $R$-submodules of $V$ w.r.t. $V[\frak{T}_{ls}^{r}(W)]$ in case $%
R_{R}$ is an injective cogenerator\ ($R$ a QF ring).

\begin{theorem}
\label{lrs-bet}Let $P=(V,W)$ be a left $R$-pairing and consider $V$ with the
linear weak topology $V[\frak{T}_{ls}^{r}(W)].$ Assume $R_{R}$ to be an
injective cogenerator.

\begin{enumerate}
\item  The closure of a right $R$-submodule $X\subseteq V$ is given by $%
\overline{X}=X^{\bot \bot }.$

\item  Let $X\subset Y\subseteq V$ be right $R$-submodules. Then $X$ is
dense in $Y$ if and only if $X^{\bot }=Y^{\bot }.$ If $W\overset{\chi _{P}}{%
\hookrightarrow }V^{\ast },$ then $X\subset V$ is dense if and only if $%
X^{\bot }=0.$

\item  Let $R$ be a QF ring and $X\subset V$ be an $R$-cofinite right $R$%
-submodule. Then $X$ is closed if and only if $\mathrm{An}(X)=\chi
_{P}(X^{\bot }).$

\item  The class of closed $R$-submodules of $V$ is given by 
\begin{equation*}
{\normalsize \{K^{\bot }|}\text{ }K\subset W\text{ is an \emph{arbitrary}
left }R\text{-submodule}\}.
\end{equation*}

\item  If $R$ is a QF-ring and $W\overset{\chi _{P}}{\hookrightarrow }%
V^{\ast }$ is an embedding, then the class of open $R$-submodules of $V$ is
given by 
\begin{equation*}
{\normalsize \{K^{\bot }|}\text{ }K\subset W\text{ is a \emph{finitely
generated} left }R\text{-submodule}\}.
\end{equation*}
\end{enumerate}
\end{theorem}

\begin{Beweis}
\begin{enumerate}
\item  By Lemma \ref{orth-clos} (1) $\overline{X}\subseteq X^{\bot \bot }.$
On the other hand, let $\widetilde{v}\in X^{\bot \bot }\backslash \overline{X%
}$ be arbitrary. Then there exists by \ref{An-Ke} (1) a finitely generated
left $R$-submodule $K\subset W,$ such that $\widetilde{v}\notin X+K^{\bot }=%
\mathrm{KeAn}(X+K^{\bot }).$ Consequently there exists $\delta \in V^{\ast },
$ such that $\delta (X+K^{\bot })=0$ and $\delta (\widetilde{v})\neq 0$. By
assumption $R_{R}$ is injective and it follows from \ref{An-Ke} (3) that $%
\delta \in \mathrm{An}(K^{\bot })=\mathrm{AnKe}(\chi _{P}(K))=\chi _{P}(K),$
i.e. $\delta =\chi _{P}(w)$ for some $w\in K.$ So 
\begin{equation*}
0=<\widetilde{v},w>=\chi _{P}(w)(\widetilde{v})=\delta (\widetilde{v})\neq 0,
\end{equation*}
a contradiction. It follows then that $\overline{X}=X^{\bot \bot }.$

\item  $X\subset Y$ is dense if and only if $\overline{X}=\overline{Y}$ and
the result follows from (1).

\item  Let $R$ be a QF ring and $X\subseteq V$ be an $R$-cofinite right $R$%
-submodule. Let $X$ be closed, i.e. $X=X^{\bot \bot }$ by (1). Since $_{R}R$
is Noetherian, $\chi _{P}(X^{\bot })\subseteq \mathrm{An}(X)\simeq
(V/X)^{\ast }$ is finitely generated in $_{R}\mathcal{M}.$ Since $R_{R}$ is
injective, we have by \ref{An-Ke} (3): 
\begin{equation*}
\mathrm{An}(X)=\mathrm{AnKeAn}(X)=\mathrm{AnKeAn}(X^{\perp \perp })=\mathrm{%
AnKeAn}(\mathrm{Ke}(\chi _{P}(X^{\perp })))=\chi _{P}(X^{\perp }).
\end{equation*}
On the other hand, if $\mathrm{An}(X)=\chi _{P}(X^{\perp }),$ then it
follows by Lemma \ref{orth-clos} (3) that $X$ is closed and we are done.

\item  Follows from (1) and Lemma \ref{orth-clos} (1).

\item  Let $R$ be a QF ring and $W\overset{\chi _{P}}{\hookrightarrow }%
V^{\ast }.$ If $K\subset W$ is a finitely generated left $R$-submodule, then 
$K^{\bot }\subset V$ is open by definition. On the other hand, if $X\subset V
$ is an open right $R$-submodule, then $X$ is closed, i.e. $X=X^{\bot \bot }$%
. By Lemma \ref{orth-clos} (2) $X\subset V$ is $R$-cofinite and so $X^{\bot }%
\overset{\chi _{P}}{\hookrightarrow }\mathrm{An}(X)\simeq (V/X)^{\ast }$ is
finitely generated in $_{R}\mathcal{M}.\blacksquare $
\end{enumerate}
\end{Beweis}

\begin{corollary}
\label{th-stet}Let $(V,W),(V^{\prime },W^{\prime })$ be left $R$-pairings
and consider $V$ and $V^{\prime }$ with the linear weak topology $V[\frak{T}%
_{ls}^{r}(W)],$ $V^{\prime }[\frak{T}_{ls}^{r}(W^{\prime })]$ respectively.
Let $(\xi ,\theta ):(V^{\prime },W^{\prime })\rightarrow (V,W)$ be a
morphism of left $R$-pairings.

\begin{enumerate}
\item  If $K^{\prime }\subset W^{\prime }$ is a left $R$-submodule, then $%
\xi ^{-1}(K^{\prime \bot })=(\theta (K^{\prime }))^{\bot }.$ In particular, $%
\xi :V\rightarrow V^{\prime }$ is continuous. In particular, $\xi
^{-1}(Y^{\prime })\subset V$ is closed for every closed right $R$-submodule $%
Y^{\prime }\subset V^{\prime }.$

\item  If $R_{R}$ is an injective cogenerator, then $\xi ^{-1}(Y^{\prime
})\subset V$ is orthogonally closed for every closed right $R$-submodule $%
Y^{\prime }\subset V^{\prime }.$
\end{enumerate}
\end{corollary}

\begin{Beweis}
\begin{enumerate}
\item  Trivial.

\item  If $Y^{\prime }\subset V^{\prime }$ is closed, then it follows by
Theorem \ref{lrs-bet} (3) that $Y^{\prime }=K^{\prime \bot }$ for some $R$%
-submodule $K^{\prime }\subset W^{\prime }.$ It follows then by (1) that $%
\xi ^{-1}(Y^{\prime })=\xi ^{-1}(K^{\prime \bot })=(\theta (K^{\prime
}))^{\bot },$ i.e. $\xi ^{-1}(Y^{\prime })\subset V$ is orthogonally closed.$%
\blacksquare $
\end{enumerate}
\end{Beweis}

\begin{proposition}
\label{f*-clos}Let $W,W^{\prime }$ be left $R$-modules and consider $^{\ast
}W,$ $^{\ast }W^{\prime }$ with the finite topology. Let $\theta \in \mathrm{%
Hom}_{R-}(W^{\prime },W)$ and consider the morphism of left $R$-pairings $%
(\theta ^{\ast },\theta ):(W^{\prime \ast },W^{\prime })\rightarrow (W^{\ast
},W).$

\begin{enumerate}
\item  $\theta ^{\ast -1}(\mathrm{An}(K^{\prime }))=\mathrm{An}(\theta
(K^{\prime }))$ for every left $R$-submodule $K^{\prime }\subset W^{\prime }.
$ In particular $\theta ^{\ast }:$ $^{\ast }W\rightarrow $ $^{\ast
}W^{\prime }$ is continuous.

\item  If $_{R}R$ is $W$-injective, then $\theta ^{\ast }(\mathrm{An}(K))=%
\mathrm{An}(\theta ^{-1}(K))$ for every left $R$-submodule $K\subseteq W.$

\item  If $R_{R}$ is an injective cogenerator and $_{R}R$ is $W$-injective 
\emph{(}e.g. $R$ is a QF-ring\emph{)}, then

\begin{enumerate}
\item  $\theta ^{\ast }:$ $^{\ast }W\rightarrow $ $^{\ast }W^{\prime }$ is 
\emph{linearly closed} \emph{(}i.e. $\theta ^{\ast }(X)\subset $ $^{\ast
}W^{\prime }$ is closed for every closed right $R$-submodule $X\subset $ $%
^{\ast }W$\emph{)}.

\item  $\overline{\theta ^{\ast }(X)}=\theta ^{\ast }(\overline{X})$ for
every right $R$-submodule $X\subset $ $^{\ast }W.$

\item  $\mathrm{Ke}(\theta ^{\ast }(X))=\theta ^{-1}(\mathrm{Ke}(X))$ for
every right $R$-submodule $X\subset $ $^{\ast }W.$

\item  For $R$-submodules $X_{1},...,X_{k}\subset $ $^{\ast }W$ we have $%
\overline{X_{1}+...+X_{k}}=\overline{X_{1}}+...+\overline{X_{k}}.$ Hence
every \emph{finite} sum of closed right $R$-submodules of $^{\ast }W$ is
closed.
\end{enumerate}
\end{enumerate}
\end{proposition}

\begin{Beweis}
\begin{enumerate}
\item  Trivial.

\item  Let $K\subseteq W$ be a left $R$-submodule. Clearly $\theta ^{\ast }(%
\mathrm{An}(K))\subseteq \mathrm{An}(\theta ^{-1}(K)).$ On the other hand,
consider the $R$-linear mapping 
\begin{equation}
0\rightarrow W^{\prime }/\theta ^{-1}(K)\overset{\iota }{\hookrightarrow }%
W/K.  \label{th-1}
\end{equation}
By assumption $_{R}R$ is $W$-injective and so it is $W/K$-injective (e.g. 
\cite[16.2]{Wis91}). Hence (\ref{th-1}) induces the epimorphism 
\begin{equation*}
^{\ast }(W/K)\overset{\iota ^{\ast }}{\longrightarrow }\text{ }^{\ast
}(W^{\prime }/\theta ^{-1}(K))\longrightarrow 0,
\end{equation*}
or equivalently the epimorphism 
\begin{equation*}
\mathrm{An}(K)\overset{\theta ^{\ast }}{\longrightarrow }\mathrm{An}(\theta
^{-1}(K))\longrightarrow 0.
\end{equation*}

\item  Let $R_{R}$ be an injective cogenerator and $_{R}R$ be $W$-injective.

\begin{enumerate}
\item  The result follows from Theorem \ref{lrs-bet} (1), Lemma \ref
{orth-clos} (1) and (2).

\item  Let $X\subset $ $^{\ast }W$ be a right $R$-submodule. By (a) $\theta
^{\ast }$ is linearly closed, so $\overline{\theta ^{\ast }(X)}\subseteq
\theta ^{\ast }(\overline{X}).$ By (1) $\theta ^{\ast -1}(\overline{\theta
^{\ast }(X)})$ is closed and it follows that $\overline{X}\subseteq \theta
^{\ast -1}(\overline{\theta ^{\ast }(X)}),$ i.e. $\theta ^{\ast }(\overline{X%
})\subseteq \overline{\theta ^{\ast }(X)}$ and the result follows.

\item  For every right $R$-submodule $X\subset $ $^{\ast }W$ we get by the
results above: 
\begin{equation*}
\begin{tabular}{lllll}
$\mathrm{Ke}(\theta ^{\ast }(X))$ & $=$ & $\mathrm{KeAnKe}(\theta ^{\ast
}(X))$ & $=$ & \textrm{Ke}$(\overline{\theta ^{\ast }(X)})$ \\ 
& $=$ & \textrm{Ke}$(\theta ^{\ast }(\overline{X}))$ & $=$ & \textrm{Ke}$%
(\theta ^{\ast }(\mathrm{AnKe}(X)))$ \\ 
& $=$ & $\theta ^{-1}(\mathrm{KeAnKe}(X))$ & $=$ & $\theta ^{-1}(\mathrm{Ke}%
(X)).$%
\end{tabular}
\end{equation*}

\item  Let $X_{1},...,X_{k}\subset $ $^{\ast }W$ be right $R$-submodules. By
Theorem \ref{lrs-bet} (1) and induction on $k$ in \ref{An-Ke} (2) we have 
\begin{equation*}
\overline{\sum_{i=1}^{k}X_{i}}=\mathrm{AnKe}(\sum_{i=1}^{k}X_{i})=\mathrm{An}%
(\bigcap\limits_{i=1}^{k}\mathrm{Ke}(X_{i}))=\sum_{i=1}^{k}\mathrm{AnKe}%
(X_{i})=\sum_{i=1}^{k}\overline{X_{i}}.\blacksquare 
\end{equation*}
\end{enumerate}
\end{enumerate}
\end{Beweis}

\section{The $\protect\alpha$-condition}

\qquad In a joint work with J. G\'{o}mez-Torrecillas and J. Lobillo \cite
{AG-TL2001} we presented the so called $\alpha $\emph{-condition }for
pairings\emph{\ }over commutative rings, which has shown to be a natural
assumption in the author's study of \emph{duality theorems }for Hopf
algebras \cite{Abu2001}. Recently that condition has shown to be a natural
assumption in the study of the category of right (left) comodules of a \emph{%
coring} $\mathcal{C}$ as a full subcategory of the category right (left)
modules of its \emph{dual ring} $^{\ast }\mathcal{C}$ (e.g. \cite{Abu03}).
In this section we consider this condition for pairings over arbitrary (not
necessarily commutative) rings and give examples of pairings satisfying it.
In particular we extend our observations in \cite{Abu2001} on such pairings
from the commutative case to the arbitrary one.

\begin{punto}
\textbf{The category }$\mathcal{P}_{l}^{\alpha }.$\label{P-alph} We say a
left $R$-pairing $P=(V,W)$ satisfies the $\alpha $\emph{-condition}\textbf{\ 
}(or $P$ is an $\alpha $\emph{-pairing}) iff for every right $R$-module $M$
the following mapping is injective 
\begin{equation}
\alpha _{M}^{P}:M\otimes _{R}W\rightarrow \mathrm{Hom}_{-R}(V,M),\text{ }%
\sum m_{i}\otimes w_{i}\mapsto \lbrack v\mapsto \sum m_{i}<v,w_{i}>].
\label{alp}
\end{equation}

With $\mathcal{P}_{l}^{\alpha }\subset \mathcal{P}_{l}$ we denote the \emph{%
full} subcategory of left $R$-pairings satisfying the $\alpha $-condition
(we call these \emph{left }$\alpha $\emph{-pairings}). We call a left $R$%
-pairing $P=(V,W)$ \emph{dense}, if $\kappa _{P}(V)\subseteq $ $^{\ast }W$
is dense w.r.t. the finite topology. The subcategory of \emph{right }$\alpha 
$\emph{-pairings} $\mathcal{P}_{r}^{\alpha }\subset \mathcal{P}_{r}$ is
defined analogously.
\end{punto}

\begin{remark}
\label{flat}Let $P=(V,W)\in \mathcal{P}_{l}^{\alpha }.$ Then $W\overset{\chi
_{P}}{\hookrightarrow }V^{\ast },$ hence $_{R}W$ is in particular $R$%
-cogenerated. If $M$ is an arbitrary right $R$-module, then we have for
every right $R$-submodule $N\subset M$ the commutative diagram 
\begin{equation*}
\xymatrix{ N \otimes_{R} W \ar[rr]^{\alpha _N ^P} \ar[d]_{\iota _N \otimes
id_W} & & {\rm Hom}_{-{R}} (V,N) \ar@{^{(}->}[d] \\ M \otimes_{R} W
\ar[rr]_{\alpha_M ^P} & & {\rm Hom}_{-{R}} (V,M)}
\end{equation*}
By assumption $\alpha _{N}^{P}$ is injective and so $N\subset M$ is $W$%
-pure. By our choice $M$ is an arbitrary $R$-module, hence $_{R}W$ is flat.
If $_{R}W$ is finitely presented or $R$ is left perfect, then $_{R}W$ is
projective.
\end{remark}

An important observation for $\alpha $-pairings is

\begin{lemma}
\label{q-2}Let $P=(V,W)\in \mathcal{P}_{l}^{\alpha }.$ For every right $R$%
-module $M$ and every $R$-submodule $N\subset M$ we have for arbitrary $\sum
m_{i}\otimes w_{i}\in M\otimes _{R}W:$%
\begin{equation}
\sum m_{i}\otimes w_{i}\in N\otimes _{R}W\Leftrightarrow \sum
m_{i}<v,w_{i}>\in N\text{ for all }v\in V.  \label{NotD}
\end{equation}
\end{lemma}

\begin{Beweis}
By remark \ref{flat}, $_{R}W$ is flat and we get the commutative diagram
with exact rows 
\begin{equation*}
\xymatrix{0 \ar[r] & N \otimes_R W \ar[rr]^{\iota_N \otimes_R id_W}
\ar@{^{(}->}[d]^{\alpha_N ^P} & & M \otimes_R W \ar[rr]^{\pi \otimes id_W}
\ar@{^{(}->}[d]^{\alpha_M ^P} & & M/N \otimes_R W \ar@{^{(}->}[d]^{\alpha
_{M/N} ^P} \ar[r] & 0\\ 0 \ar[r] & {\rm Hom} _{-R}(V,N)
\ar[rr]_{(V,\iota_N)} & & {\rm Hom} _{-R}(V,M) \ar[rr]_{(V,\pi)} & & {\rm
Hom} _{-R}(V,M/N) & }
\end{equation*}
Obviously $\sum m_{i}<v,w_{i}>\in N$ for all $v\in V$ if and only if 
\begin{equation*}
\sum m_{i}\otimes w_{i}\in \mathrm{Ker}{\normalsize ((V,\pi )\circ \alpha
_{M}^{P})=\mathrm{Ker}(\alpha _{M/N}^{P}\circ (\pi \otimes id_{W}))=\mathrm{%
Ker}(\pi \otimes id_{W})=N\otimes _{R}W.\blacksquare }
\end{equation*}
\end{Beweis}

\qquad

\begin{proposition}
\label{rp-rp}

\begin{enumerate}
\item  Let $P=(V,W)$ be a left $R$-pairing.

\begin{enumerate}
\item  Let $W^{\prime }\subset W$ be a left $R$-submodule and consider the
induced left $R$-pairing $P^{\prime }:=(V,W^{\prime }).$ If $P^{\prime }\in 
\mathcal{P}_{l}^{\alpha },$ then $W^{\prime }\subset W$ is pure. If $P\in 
\mathcal{P}_{l}^{\alpha },$ then $P^{\prime }\in \mathcal{P}_{l}^{\alpha }$
if and only if $W^{\prime }\subset W$ is pure.

\item  Let $V^{\prime }\subset V$ be a right $R$-submodule, $W^{\prime
}\subset W$ be a left $R$-submodule with $<V^{\prime },W^{\prime }>=0$ and
consider the left $R$-subpairing $Q:=(V/V^{\prime },W^{\prime })$ of $P.$ If 
$P\in \mathcal{P}_{l}^{\alpha },$ then $Q\in \mathcal{P}_{l}^{\alpha }$ if
and only if $W^{\prime }\subset W$ is pure. In particular $\mathcal{P}%
_{l}^{\alpha }$ is closed under pure left $R$-subpairings.
\end{enumerate}

\item  Let $\Omega =(Y,W)$ be a left $R$-pairing, $V$ be a right $R$-module, 
$\xi :V\rightarrow Y$ be an $R$-linear mapping, $P:=(V,W)$ be the induced
left $R$-pairing and consider the following statements:

(i) $\Omega \in \mathcal{P}_{l}^{\alpha }$ and $P$ is dense;

(ii) $\Omega \in \mathcal{P}_{l}^{\alpha }$ and $\xi (V)\subset Y$ is dense
w.r.t. $Y[\frak{T}_{_{ls}}^{r}(W)];$

(iii) $P\in \mathcal{P}_{l}^{\alpha };$

(iv) $P\in \mathcal{P}_{l}^{\alpha }$ and $W\overset{\chi _{P}}{%
\hookrightarrow }V^{\ast }$ is an embedding.

The following implications are always true: (i) $\Longrightarrow $ (ii) $%
\Longrightarrow $ (iii) $\Longrightarrow $ (iv). If $R_{R}$ is an injective
cogenerator, then (i)-(iv)\ are equivalent.
\end{enumerate}
\end{proposition}

\begin{Beweis}
\begin{enumerate}
\item  The result follows from the commutativity of the following diagram
for every right $R$-module $M$%
\begin{equation*}
\xymatrix{ M \otimes_R W' \ar[d]_{id_M \otimes {\iota}_{W'}}
\ar[drr]_(.4){\alpha_M ^{P'}} \ar[rr]^(.45){\alpha_M ^{Q}} & & {\rm
Hom}_{-R} (V/V',M) \ar@{^{(}->}[d]\\ M \otimes_{R} W \ar[rr]_{\alpha_M ^{P}}
& & {\rm Hom} _{-R} (V,M)}
\end{equation*}

\item  Consider for every right $R$-module $M$ the commutative diagram 
\begin{equation*}
\xymatrix{ M \otimes_R W \ar[rr]^{\alpha _M ^{\Omega}} \ar[drr]_{\alpha_{M}
^{P}} & & {\rm Hom} _{-R}(Y,M) \ar[d]^{(\xi ,M)} \\ & & {\rm Hom} _{-R}
(V,M)}
\end{equation*}
(i) $\Rightarrow $ (ii) trivial.

(ii) $\Rightarrow $ (iii) Let $\Omega \in \mathcal{P}_{l}^{\alpha }$ and
assume that $\xi (V)\subset Y$ is dense. Let $\sum\limits_{i=1}^{n}m_{i}%
\otimes w_{i}\in \mathrm{Ker}(\alpha _{M}^{P}).$ By assumption for every $%
y\in Y$ there exists some $v_{y}\in V,$ such that $\kappa _{\Omega
}(y)(w_{i})=\kappa _{P}(v_{y})(w_{i})$ for $i=1,...,n$ and it follows then
that 
\begin{equation*}
\begin{tabular}{lllll}
$\alpha _{M}^{\Omega }(\sum\limits_{i=1}^{n}m_{i}\otimes w_{i})(y)$ & $=$ & $%
\sum\limits_{i=1}^{n}m_{i}<y,w_{i}>$ & $=$ & $\sum%
\limits_{i=1}^{n}m_{i}<v_{y},w_{i}>$ \\ 
& $=$ & $\alpha _{M}^{P}(\sum\limits_{i=1}^{n}m_{i}\otimes w_{i})(v_{y})$ & $%
=$ & $0.$%
\end{tabular}
\end{equation*}
So $\mathrm{Ker}(\alpha _{M}^{P})=\mathrm{Ker}(\alpha _{M}^{\Omega })=0,$
i.e. $\alpha _{M}^{P}$ injective. The $R$-module $M$ is by our choice
arbitrary and so $P\in \mathcal{P}_{l}^{\alpha }.$

(iii) $\Rightarrow $ (iv) Trivial.

Let $R_{R}$ be an injective cogenerator.

(iv) $\Longrightarrow $ (i) If $W\overset{\chi _{P}}{\hookrightarrow }%
V^{\ast }$ is an embedding, then it follows by Theorem \ref{lrs-bet} (1)
that $\overline{\kappa _{P}(V)}=\mathrm{AnKe}(\kappa _{P}(V))=\mathrm{An}%
(V^{\bot })=\mathrm{An}\{0_{W}\}=$ $^{\ast }W,$ i.e. $P$ is a dense left $R$%
-pairing.$\blacksquare $
\end{enumerate}
\end{Beweis}

Over Noetherian rings we have the following interesting observation:

\begin{proposition}
\label{P-rs}Let $V$ be a right $R$-module, $R[V]$ be the free right $R$%
-module with basis $V,$ $W\subset V^{\ast }$ be a left $R$-submodule and
consider the left $R$-pairing $P:=(V,W).$ Assume $R_{R}$ to be Noetherian.

\begin{enumerate}
\item  For every right $R$-module $M$ the following mapping is injective 
\begin{equation}
\beta _{M}:M\otimes _{R}R^{V}\rightarrow M^{V},\text{ }m\otimes f\mapsto
\lbrack v\mapsto mf(v)],  \label{bet_M}
\end{equation}
i.e. $\widetilde{P}:=(R[V],R^{V})$ is a left $\alpha $-pairing.

\item  Let $M$ be an arbitrary right $R$-module. Then the canonical mapping $%
\alpha _{M}^{P}:M\otimes _{R}W\rightarrow \mathrm{Hom}_{-R}(V,M)$ is
injective if and only if $W\subset R^{V}$ is $M$-pure. If moreover $V_{R}$
is projective, then $\alpha _{M}^{P}$ is injective if and only if $%
W\subseteq V^{\ast }$ is $M$-pure.

\item  The following statements are equivalent:

(i) $P\in \mathcal{P}_{l}^{\alpha };$

(ii) $\alpha _{M}^{P}$ is injective for every \emph{(}finitely presented%
\emph{)} right $R$-module $M;$

(iii) $W\subset R^{V}$ is pure.
\end{enumerate}
\end{proposition}

\begin{Beweis}
\begin{enumerate}
\item  Let $M$ be an arbitrary right $R$-module and write $M$ as a direct
limit of its finitely generated $R$-submodules $M={\underrightarrow{lim}}%
_{\Lambda }M_{\lambda }$ (e.g. \cite[24.7]{Wis91}). For every $\lambda \in
\Lambda $ the $R$-module $M_{\lambda }$ is finitely presented in $\mathcal{M}%
_{R}\ $and so 
\begin{equation*}
\beta _{M_{\lambda }}:M_{\lambda }\otimes _{R}R^{V}\rightarrow M_{\lambda
}^{V}
\end{equation*}
is an isomorphism (e.g. \cite[25.4]{Wis91}). Moreover, for every $\lambda
\in \Lambda $ the restriction of $\beta _{M}$ on $M_{\lambda }$ coincides
with $\beta _{M_{\lambda }}$ and so the following mapping is injective: 
\begin{equation*}
\beta _{M}=\underrightarrow{lim}\beta _{M_{\lambda }}:\underrightarrow{lim}%
M_{\lambda }\otimes _{R}R^{V}\rightarrow \underrightarrow{lim}M_{\lambda
}^{V}\subset M^{V}
\end{equation*}
Obviously $\widetilde{P}\in \mathcal{P}_{l}^{\alpha }$ if and only if $\beta
_{M}$ is injective is for every $M\in \mathcal{M}_{R}.$

\item  The first statement follows by (1). If moreover $V_{R}$ is
projective, then the exact sequence $R[V]\rightarrow V\rightarrow 0$ splits,
hence $V^{\ast }\subset R^{V}$ is pure (direct summand) and we are done.

\item  By \cite[34.5]{Wis91}, $W\subset R^{V}$ is pure if and only if $%
W\subset R^{V}$ is $M$-pure for every finitely presented right $R$-module $M.
$ The result follows then from (2).$\blacksquare $
\end{enumerate}
\end{Beweis}

\begin{definition}
The ring $R$ is called right (\emph{semi})\emph{\ hereditary} iff every
(finitely generated) right ideal is projective.
\end{definition}

\begin{lemma}
\label{hered}Let $R_{R}$ be Noetherian and hereditary and $V$ be a right $R$%
-module. Then:

\begin{enumerate}
\item  $\breve{P}:=(V,V^{\ast })\in \mathcal{P}_{l}^{\alpha }.$

\item  Let $W\subseteq V^{\ast }$ be a left $R$-module and $P:=(W,V).$ Then $%
P\in \mathcal{P}_{l}^{\alpha }$ if and only if $W\subset V^{\ast }$ is a
pure $R$-submodule.
\end{enumerate}
\end{lemma}

\begin{Beweis}
Assume $R_{R}$ to be Noetherian and hereditary. It follows then by 
\cite[26.6]{Wis91} that $_{R}R^{\Lambda }$ is flat for every set $\Lambda .$
Moreover we have by \cite[39.13]{Wis91} that all $R$-cogenerated left $R$%
-modules are flat. Consider now $R[V],$ the free right $R$-module with basis 
$V,$ and the exact sequence of right $R$-modules 
\begin{equation}
0\longrightarrow \mathrm{Ker}(\pi )\overset{\iota }{\longrightarrow }R[V]%
\overset{\pi }{\longrightarrow }V\longrightarrow 0,  \label{Ke(pi)}
\end{equation}
with $\iota $ the embedding map and $\pi $ the canonical epimorphism. Then (%
\ref{Ke(pi)}) induces the exact sequence of left $R$-modules 
\begin{equation}
0\longrightarrow V^{\ast }\overset{\pi ^{\ast }}{\longrightarrow }R^{V}%
\overset{\iota ^{\ast }}{\longrightarrow }\mathrm{\func{Im}}(\iota ^{\ast
})\longrightarrow 0.  \label{pure-seq}
\end{equation}
Since $\mathrm{\func{Im}}(\iota ^{\ast })\subseteq \mathrm{Ker}(\pi )^{\ast
},$ $\mathrm{\func{Im}}(\iota ^{\ast })$ is an $R$-cogenerated left $R$%
-module, hence flat. Consequently $V^{\ast }\hookrightarrow R^{V}$ is pure
(e.g. \cite[36.6]{Wis91}). By Proposition \ref{P-rs} (1) the canonical
mapping $\beta _{M}:M\otimes _{R}R^{V}\rightarrow M^{V}$ is injective for
every $M\in \mathcal{M}_{R}$ and the result follows then from the
commutativity of the following diagram 
\begin{equation*}
\xymatrix{ M \otimes_R W \ar[rr]^(.45){\alpha_M ^P} \ar[d]_{id \otimes
\iota_W} & & {\rm Hom} _{-R} (V,M) \ar@{^{(}->}[r] & M^V\\ M \otimes_R V^*
\ar[urr]^(.4){\alpha_M ^{\breve{P}}} \ar@{^{(}->}[rr] & & M \otimes_R R^V
\ar@{^{(}->}[ur]_{\beta_M} & }
\end{equation*}
\end{Beweis}

\begin{lemma}
\label{p-2}Let $V,W$ be $R$-bimodules.

\begin{enumerate}
\item  If $P=(V,W),$ $P^{\prime }=(V^{\prime },W^{\prime })$ are left $%
\alpha $-pairings, then $P\otimes _{l}P^{\prime }:=(V^{\prime }\otimes
_{R}V,W\otimes _{R}W^{\prime })$ is a left $\alpha $-pairing with 
\begin{equation*}
\kappa _{P\otimes _{l}P^{\prime }}(v^{\prime }\otimes v)(w\otimes w^{\prime
})=<v,w<v^{\prime },w^{\prime }>>=<<v^{\prime },w^{\prime }>v,w>.
\end{equation*}

\item  If $P=(V,W),$ $P^{\prime }=(V^{\prime },W^{\prime })$ are right $%
\alpha $-pairings, then $P\otimes _{r}P^{\prime }:=(V\otimes _{R}V^{\prime
},W^{\prime }\otimes _{R}W)$ is a right $\alpha $-pairing with 
\begin{equation*}
\kappa _{P^{\prime }\otimes _{r}P}(v\otimes v^{\prime })(w^{\prime }\otimes
w)=<v,<v^{\prime },w^{\prime }>w>=<v<v^{\prime },w^{\prime }>,w>.
\end{equation*}
\end{enumerate}
\end{lemma}

\begin{Beweis}
We prove (1). The proof of (2) is similar. For arbitrary $M\in \mathcal{M}%
_{R}$ consider the following commutative diagram 
\begin{equation*}
\xymatrix{ M \otimes_R W \otimes_R W' \ar[rr]^(.45){\alpha_M ^{P \otimes
P'}} \ar[d]_{\alpha_{M \otimes_R W} ^{P'}} & & {\rm Hom}_{-R} (V \otimes_R
V',M) \\ {\rm Hom}_{-R} (V',M \otimes_R W) \ar[rr]_(.45){(V',\alpha_M ^P)} &
& {\rm Hom}_{-R} (V',{\rm Hom}_{-R} (V,M)) \ar[u]_{\zeta ^l}}
\end{equation*}
where $\zeta ^{l}$ is the canonical isomorphism. By assumption the $R$%
-linear mappings $\alpha _{M\otimes _{R}W}^{P^{\prime }}$ and $\alpha
_{M}^{P}$ are injective and so $\alpha _{M}^{P\otimes P^{\prime }}$ is
injective. The last statement is obvious.$\blacksquare $
\end{Beweis}

\begin{corollary}
\label{uno}Let $R_{R}$ be Noetherian.

\begin{enumerate}
\item  Let $X,X^{\prime }$ be sets, $E\subseteq R^{X}$ be a right $R$%
-submodule and $E^{\prime }\subseteq R^{X^{\prime }}$ be a left $R$%
-submodule. If $E^{\prime }\subseteq R^{X^{\prime }}$ is $E$-pure, then the
following mapping is injective: 
\begin{equation}
\delta :E\otimes _{R}E^{\prime }\rightarrow (R^{X})^{X^{\prime }},\text{ }%
f\otimes f^{\prime }\mapsto \lbrack (x,x^{\prime })\mapsto f(x)f^{\prime
}(x^{\prime })].  \label{pi}
\end{equation}

\item  Let $W,W^{\prime }$ be $R$-bimodules, $X\subset $ $^{\ast
}W,X^{\prime }\subset $ $^{\ast }W^{\prime }$ be $R$-subbimodules and
consider the canonical $R$-linear mappings 
\begin{equation*}
\kappa :X^{\prime }\otimes _{R}X\rightarrow \text{ }^{\ast }(W\otimes
_{R}W^{\prime })\text{ and }\chi :W\otimes _{R}W^{\prime }\rightarrow
(X^{\prime }\otimes _{R}X)^{\ast }.
\end{equation*}
If $W_{R}$ is flat and $\mathrm{Ke}(X)_{R}\subset $ $W_{R}$ is pure, then 
\begin{equation}
\mathrm{Ke}(\kappa (X^{\prime }\otimes _{R}X))\simeq \mathrm{Ke}(X)\otimes
_{R}W^{\prime }+W\otimes _{R}\mathrm{Ke}({X^{\prime }}).  \label{KeX-Y}
\end{equation}
\newline
\end{enumerate}
\end{corollary}

\begin{Beweis}
\begin{enumerate}
\item  Since $R_{R}$ is coherent, $_{R}R^{X^{\prime }}$ is flat in $_{R}%
\mathcal{M}$ by (\cite[26.6]{Wis91}). The result follows then from
Proposition \ref{P-rs} (1).

\item  Consider the embeddings $E:=W/\mathrm{Ke}(X)\hookrightarrow X^{\ast },
$ $E^{\prime }:=W^{\prime }/\mathrm{Ke}({X^{\prime })}\hookrightarrow
R^{X^{\prime }}$ and the commutative diagram 
\begin{equation*}
\xymatrix{ W \otimes_{R} W' \ar[rrr]^{\chi } \ar[dd]_{\pi \otimes \pi '} & &
& (X' \otimes_{R} X)^* \ar@{^{(}->}[dd]^{\iota} \\ & W/{\rm
{Ke}(X)}\otimes_{R} {R}^{X'} \ar@{^{(}->}[r] & X^* \otimes_{R} {R}^{X'}
\ar@{^{(}->}[rd]^{\beta_{X^*}} & \\ W/{\rm {Ke}(X)}\otimes_{R}
W'/{\rm{Ke}(X')} \ar@{^{(}->}[ur] \ar@{.>}[rrr]_(.6){\delta} & & &
({X^*})^{X'} }
\end{equation*}
It follows by assumptions that $W/\mathrm{Ke}(X)$ is flat in $\mathcal{M}_{R}
$ and $_{R}R^{X^{\prime }}$ is flat (e.g. \cite[36.5, 26.6]{Wis91}).
Moreover $\beta _{X^{\ast }}$ is injective by Lemma \ref{P-rs}, hence $%
\delta $ is injective. It follows then by \cite[II-3.6]{Bou74} that 
\begin{equation*}
\begin{tabular}{lllll}
\textrm{Ke}$(\kappa (X\otimes _{R}X^{\prime }))$ & $:=$ & $\mathrm{Ker}(\chi
)$ & $=$ & \textrm{Ker}$(\delta \circ (\pi \otimes \pi ^{\prime }))$ \\ 
& $=$ & $\mathrm{Ker}(\pi _{X}\otimes \pi _{X^{\prime }})$ & $=$ & \textrm{Ke%
}$(X)\otimes _{R}W^{\prime }+W\otimes _{R}\mathrm{Ke}({X^{\prime }}%
).\blacksquare $%
\end{tabular}
\end{equation*}
\end{enumerate}
\end{Beweis}

\begin{punto}
\label{alph-W}We say a left (respectively a right) $R$-module $W$ satisfies
the $\alpha $\emph{-condition}, if $(^{\ast }W,W)$ (respectively $(W^{\ast
},W)$) satisfies the $\alpha $-condition. Such modules were called \emph{%
universally torsionless}\textbf{\ }by G. Garfinkel \cite{Gar76}.
\end{punto}

\begin{punto}
\textbf{Locally projective modules}. An $R$-module $W$ is called \emph{%
locally projective} (in the sense of B. Zimmermann-Huisgen \cite{Z-H76}) iff
for every diagram of $R$-modules 
\begin{equation*}
\xymatrix{0 \ar[r] & F \ar@{.>}[dr]_{g' \circ \iota} \ar[r]^{\iota} & W
\ar[dr]^{g} \ar@{.>}[d]^{g'} & & \\ & & L \ar[r]_{\pi} & N \ar[r] & 0}
\end{equation*}
with exact rows, $F$ finitely generated as an $R$-module and every $R$%
-linear mapping $g:W\rightarrow N$ there exists an $R$-linear mapping $%
g^{\prime }:W\rightarrow L,$ such that the entstanding parallelogram is
commutative. By (\cite[Theorem 3.2]{Gar76}, \cite[Theorem 2.1]{Z-H76}), $W$
is locally projective if and only if $W$ satisfies the $\alpha $-condition.
It follows directly from the definition that every projective $R$-module is
locally projective, hence satisfies the $\alpha $-condition.
\end{punto}

Before proceeding, we would like to remark that some of following results on
the $\alpha $-condition and locally projective modules appeared in the
recent manuscript \cite[42.9-42.12]{BW03}.

\begin{proposition}
\label{proj-gut}Let $W$ be a left $R$-module.

\begin{enumerate}
\item  If $_{R}W$ is locally projective, then every pure left $R$-submodule $%
K\subset W$ is locally projective. If $_{R}R$ is $W$-injective, then every
locally projective $R$-submodule of $W$ is a pure left $R$-submodule.

\item  Let $R_{R}$ be Noetherian. Then $_{R}W$ is locally projective if and
only if $_{R}W\subset $ $_{R}R^{^{\ast }W}$ is a pure $R$-submodule.
\end{enumerate}
\end{proposition}

\begin{Beweis}
\begin{enumerate}
\item  Standard.

\item  This follows from Propositions \ref{P-rs} (2).$\blacksquare $
\end{enumerate}
\end{Beweis}

\begin{corollary}
\label{PW}If $R_{R}$ is an injective cogenerator, then for every left $R$%
-pairing $(V,W)$ the following statements are equivalent:

(i) $_{R}W$ is locally projective and $P$ is dense.

(ii) $W$ satisfies the $\alpha $-condition and $P$ is dense.

(ii) $(V,W)$ is a left $\alpha $-pairing.

(iii) $W$ satisfies the $\alpha $-condition and $W\hookrightarrow V^{\ast }.$

If $R$ a QF ring, then (i)-(iii) are moreover equivalent to:

(iv) $_{R}W$ is projective and $W\hookrightarrow V^{\ast }.$

(v) $W\subset R^{V}$ is a pure $R$-submodule.
\end{corollary}

\qquad Over semisimple rings we recover the characterizations of dense
pairings over commutative base fields:

\begin{corollary}
\label{dicht=alp}Let $P=(V,W)$ be a left $R$-pairing. If $R$ is semisimple,
then 
\begin{equation*}
P\text{ is dense }\Leftrightarrow \text{ }W\subset V^{\ast }\Leftrightarrow P%
\text{ is a left }\alpha \text{-pairing.}
\end{equation*}
\end{corollary}

\textbf{Acknowledgment.} The author would like to thank the referee for his
or her remarks and for pointing out some typos in a previous version of this
note.

\end{document}